\documentclass[a4paper, 12pt]{article}

\usepackage[koi8-r,cp1251]{inputenc}
\usepackage{amsfonts,amssymb,amsmath, hyperref}
\usepackage[final]{epsfig}
\usepackage{graphicx}

\begin{document}

\begin{center}
\textbf{  A HAUSDORFF OPERATOR ON LEBESGUE SPACE WITH COMMUTING FAMILY OF PERTURBATION MATRICES IS A NON-RIESZ OPERATOR}
\end{center}

\

\begin{center}
\textbf{A. R. Mirotin}
\end{center}

\

\centerline{amirotin@yandex.ru}
%The research was supported by the grant of the State Program of Fundamental Research of Republic of Belarus

\

\textbf{Abstract} {\small We consider a generalization of Hausdorff operators on Lebesgue spaces and  under  natural conditions
prove that such an operator is not a Riesz operator provided it is non-zero. In particular, it cannot be represented as a sum of a quasinilpotent and compact operators.}

\

Key words  and phrases. Hausdorff operator,  Riesz operator, quasinilpotent operator, compact operator, Lebesgue space.

\

MSC Class: 45P05; 47G10; 47B06

\

\section{Introduction and preliminaries}

The one-dimensional Hausdorff operator
$$
(\mathcal{H}_1f)(x) =\int_\mathbb{R} f(xt)d\chi(t)
$$
 ($\chi$ is a measure  supported in $[0,1]$) was introduced by Hardy \cite[Section 11.18]{H} as a
continuous variable analog  of regular  Hausdorff transformations (or Hausdorff means) for series. Its modern $n$-dimensional generalization is as follows:
$$
(\mathcal{H}f)(x) =\int_{\mathbb{R}^m} K(u)f(A(u)x)du, \eqno(1)
$$
where $K:\mathbb{R}^m\to \mathbb{C}$ is a locally integrable function, $f:\mathbb{R}^n\to  \mathbb{C}$,  $A(u)$ stands  for a family of non-singular $n\times n$-matrices with real entries defined
on $\mathbb{R}^m$, $x\in \mathbb{R}^n$ is a column vector \cite{LL}. The modern theory of
Hausdorff operators was initiated  by Liflyand and Moricz \cite{LM}. See survey articles \cite{Ls}, \cite{CFW}  for historical remarks and the state of the art up to 2014.

Note that the map $x\mapsto Ax$ ($A\in GL(n,\mathbb{R})$) which appears in the formula (1) is the general form of an automorphism of the additive group $\mathbb{R}^n$.
This observation  leads to the  definition of a (generalized) Hausdorff  operator over a general group $G$ via the automorphisms of $G$ that was
introduced and studied by the author in \cite{JMAA}, and \cite{AddJMAA}\footnote{The  case of a Hausdorff  operator on  $p$-Adic vector spaces  was considered earlier in \cite{Vol}, the special case of a Hausdorff  operator on the Heisenberg group in the sense of this definition was considered   in \cite{RFW}, and \cite{CFWu}. The generalized Delsart translation operators (see, e.g., \cite{Lev})  lead to  Hausdorff  operators in the sense of this definition, too}. For the additive group $\mathbb{R}^n$, this definition looks as follows:

\textbf{Definition 1.} Let  $(\Omega,\mu)$ be a $\sigma$-compact  topological space endowed with a positive regular  Borel measure $\mu,$
let $K$ be a locally integrable function on $\Omega,$  and let $(A(u))_{u\in \Omega}$ be a $\mu$-measurable family of $n\times n$-matrices, non-singular
for $\mu$-almost every $u$, with $K(u) \ne 0.$
  We define the \textit{Hausdorff  operator} with the kernel $K$  by (recall that $x\in\mathbb{R}^n$ is a column vector)
$$
(\mathcal{H}_{K, A}f)(x) =\int_\Omega K(u)f(A(u)x)d\mu(u).
$$

The general form of a Hausdorff  operator given by Definition 1 (with an arbitrary measure space $(\Omega,\mu)$ instead of $\mathbb{R}^m$) gives us,
for example, the opportunity to consider (in the case $\Omega=\mathbb{Z}^m$) discrete Hausdorff  operators  (see \cite{Forum}, \cite{faa}, and Example 3 below).

The problem of compactness of  Hausdorff  operators was   posed by Liflyand \cite{L} (see also \cite{Ls}).
There is a conjecture that non-zero Hausdorff operator on $L^p(\mathbb{R}^n)$ is non-compact.
 For  the case  $p=2$ and for commuting $A(u)$ this hypothesis was confirmed in \cite{Forum} (and for the diagonal $A(u)$ --- in \cite{JMAA}).
 Moreover, we  conjectured in \cite{Conf}  that a nontrivial Hausdorff operator on $L^p(\mathbb{R}^n)$ is non-Riesz.
 %Our main result states that this conjecture is true for  a commuting family $A(u)$.

 A notion of a Riesz operator was introduced by Ruston \cite{Ruston}. Recall that a bounded operator  $T$   on some Banach space is a \textit{Riesz operator} if it possesses spectral properties like those
of a compact operator; i.~e., $T$ is a non-invertible operator whose nonzero spectrum  consists  of eigenvalues  of finite multiplicity with no limit points other then $0$.
This is equivalent to the fact that $T-\lambda$ is Fredholm for all scalars $\lambda \ne 0$ \cite{Sche}, \cite[Section 9.6]{Sche1}. For  example, a sum of a quasinilpotent and compact operator
is Riesz \cite[Theorem 3.29]{Dow}. See \cite{Sche}, \cite[Section 9.6]{Sche1}, \cite{Dow} and the bibliography therein for other interesting characterizations of Riesz operators.

In this note we prove the  aforementioned conjecture for the case where $A(u)$ is  a commuting family of   self-adjoint
 matrices. The result has been announced in \cite{arxnonR}. The case of positive or negative definite perturbation matrices was considered in \cite{Conf}.

\section{The main result}

We shall employ for lemmas to prove our main result.

\textbf{Lemma 1} \cite{JMAA} (cf. \cite[(11.18.4)]{H}, \cite{BM}). \textit{Let  $|\det A(u)|^{-1/p}K(u)\in L^1(\Omega).$
Then the operator $\mathcal{H}_{K, A}$ is bounded in   $L^p(\mathbb{R}^n)$ ($1\leq p\le \infty$) and}
$$
\|\mathcal{H}_{K, A}\|\leq \int_\Omega |K(u)||\det A(u)|^{-1/p}d\mu(u).
$$
This estimate is sharp  (see  Theorem 1  in \cite{faa}).

\textbf{Lemma 2} \cite{faa} (cf. \cite{BM}). \textit{Under the assumptions of Lemma} 1 \textit{the adjoint for the Hausdorff operator on $L^p(\mathbb{R}^n)$ is of the form}
$$
(\mathcal{H}_{K, A}^*f)(x) =\int_\Omega \overline{K(v)}|\det A(v)|^{-1}f(A(v)^{-1}x)d\mu(v).
$$
\textit{Thus, the adjoint for a Hausdorff operator  is also  Hausdorff.}

\textbf{Lemma 3.} \textit{Let   $S$ be a ball in $\mathbb{R}^n$, $q\in [1,\infty)$, and let $R_{q,S}$ denote the restriction operator
$L^q(\mathbb{R}^n)\to  L^q(S)$, $f\mapsto f|S$. If we identify the dual $(L^{q})^{*}$ of $L^q$ with $L^p$ ($1/p+1/q=1$),
then the adjoint $R_{q,S}^*$ is the operator of natural embedding $L^p(S)\hookrightarrow L^p(\mathbb{R}^n)$}.

{\it Proof.}  For $g\in L^p(S)$, let
$$
g^*(x)=\begin{cases}
g(x), &  \mathrm{for}\  x\in S,\\
0, & \mathrm{for}\  x\in \mathbb{R}^n\setminus S.
\end{cases}
$$
Then the map $g\mapsto g^*$ is the natural embedding $L^p(S)\hookrightarrow L^p(\mathbb{R}^n)$.

By definition, the adjoint  $R_{q,S}^*: L^q(S)^*\to L^q(\mathbb{R}^n)^*$ acts according to  the rule
$$
(R_{q,S}^*\Lambda)(f)=\Lambda(R_{q,S}f)\quad \ (\Lambda\in L^q(S)^*, f\in L^q(\mathbb{R}^n)).
$$
If we (by the Riesz theorem) identify the dual of $L^q(S)$ with $L^p(S)$ via the formula $\Lambda\leftrightarrow g$, where
$$
\Lambda(h)=\int_S g(x)h(x)dx\quad \ (g\in L^p(S), h\in L^q(S)),
$$
and analogously for the  dual of $L^q(\mathbb{R}^n)$, then the definition of  $R_{q,S}^*$ is of the form
$$
\int_{\mathbb{R}^n}(R_{q,S}^*g)(x)f(x)dx=\int_S g(x)(f|S)(x)dx.
$$
But
$$
\int_S g(x)(f|S)(x)dx=\int_{\mathbb{R}^n}g^*(x)f(x)dx\quad (f\in L^q(\mathbb{R}^n)).
$$
 The right-hand side of the last formula is the linear functional from $L^q(\mathbb{R}^n)^*$.
 If we (again by the Riesz theorem) identify this functional with the function $g^*$, the result follows.

 Consider the modified $n$-dimensional  Mellin transform for the $n$-hyperoctant $U$ of $\mathbb{R}^n$ in the form

$$
(\mathcal{M}f)(s):=\frac{1}{(2\pi)^{n/2}}\int_{U}|x|^{-\frac{1}{q}+is}f(x)dx,\quad s\in \mathbb{R}^n.
$$
  Here and below we assume that $1< q\le\infty$,
  $$
  x|^{-\frac{1}{q}+ is}:= \prod_{j=1}^n |x_j|^{-\frac{1}{q}+ is_j},
    $$
     where
  $$
  |x_j|^{-\frac{1}{q}+ is_j}:=\exp\left(\left(-\frac{1}{q}+ is_j\right)\log |x_j|\right).
  $$

\textbf{Lemma 4.} 1) \textit{The map $\mathcal{M}$ is a bounded operator between  $L^p(U)$ and $L^q(\mathbb{R}^n)$ for $1\leq p\leq 2$  ($1/p+1/q=1).$}

2) \textit{If we identify the dual $(L^{p})^{*}$ of $L^p$ with $L^q$ ($1/p+1/q=1$),
then the adjoint   for the operator $\mathcal{M}$ on the space $L^p(U)$ ($1\leq p\le 2$) is as follows:
 $$
(\mathcal{M}^*g)(x) =\frac{1}{(2\pi)^{n/2}}\int_{\mathbb{R}^n}|x|^{-\frac{1}{q}+is}g(s)ds,\  x\in U.
$$}

Proof. 1) It can easily be obtained from the Hausdorff--Young inequality for the $n$-dimensional  Fourier transform by using the exponential change of variables
 (see \cite{BPT}).

2) To compute $\mathcal{M}^*$, for  $g\in L^p(\mathbb{R}^n)$ consider the operator
$$
(\mathcal{M}'g)(x):=\frac{1}{(2\pi)^{n/2}}\int_{\mathbb{R}^n}|x|^{-\frac{1}{q}+is}g(s)ds,\  x\in U.
$$
This  is a bounded operator taking  $L^p(\mathbb{R}^n)$ into  $L^q(U).$ Indeed,
since
$$
|x|^{-\frac{1}{q}+ is}=\prod_{j=1}^n |x_j|^{-\frac{1}{q}}\exp(is_j\log |x_j|),
$$
we have
$$
(\mathcal{M}'g)(x)=|x|^{-\frac{1}{q}}\frac{1}{(2\pi)^{n/2}}\int_{\mathbb{R}^n}\exp(is\cdot (\log|x_j|))g(s)ds,\  x\in U,
$$
where $|x|:=|x_1|\dots |x_n|$, $(\log|x_j|):=(\log |x_1|,\dots,\log |x_n|)$, and the dot denotes the inner product on $\mathbb{R}^n$.
Thus, we can express the function  $\mathcal{M}'g$ via the Fourier transform $\widehat{g}$ of $g$ as
$$
(\mathcal{M}'g)(x)=|x|^{-1/q}\widehat{g}(-(\log|x_j|)),\quad  x\in U,
$$
and therefore
$$
\|\mathcal{M}'g\|_{L^q(U)}=\left(\int_U |x|^{-1}|\widehat{g}(-(\log|x_j|))|^q dx\right)^{1/q}.
$$
Putting in the last integral $y_j:=-\log|x_j|$ ($j=1,\dots,n$) and taking into account that the Jacobian module of this transformation is
$$
\left|\frac{\partial(x_1,\dots,x_n)}{\partial(y_1,\dots,y_n)}\right|=\det \mathrm{diag}(e^{-y_1},\dots,e^{-y_n})=\exp\left(-\sum_{j=1}^n y_j \right),
$$
we get by  the Hausdorff--Young inequality that
$$
\|\mathcal{M}'g\|_{L^q(U)}=\|\widehat{g}\|_{L^q(\mathbb{R}^n)}\le \|g\|_{L^p(\mathbb{R}^n)}.
$$

If  $f\in L^p(U)$, and  $f(x)|x|^{-1/q}\in L^1(U)$, $g\in  L^p(\mathbb{R}^n)\cap L^1(\mathbb{R}^n)$, then the Fubini--Tonelli’s theorem implies
$$
\int_{\mathbb{R}^n}(\mathcal{M}f)(s)g(s)ds=\int_{U}f(x)(\mathcal{M}'g)(x)dx.
$$
Since the bilinear  dual pairing $(\varphi,\psi)\mapsto\int \varphi \psi d\nu$ is continuous on $L^p(\nu)\times L^q(\nu)$,
the last equality is valid for all   $f\in L^p(U)$, $g\in  L^p(\mathbb{R}^n)$ by continuity. So, $\mathcal{M}'=\mathcal{M}^*$.

Now we are in position to prove our main result.

\textbf{Theorem 1.}
\textit{Let $A(u)$ be a commuting family of real self-adjoint $n\times n$-matrices  ($u$ satisfies    $K(u)\ne 0$),
and $(\det A(u))^{-1/p}K(u)\in L^1(\Omega).$ Then a Hausdorff  operator $\mathcal{H}_{K, A}$ in  $L^p(\mathbb{R}^n)$
($1\leq p\le \infty$) is  non-Riesz  (and in particular it is not a sum of a quasinilpotent and compact operator)  provided it is  non-zero.}

 Proof.
Assume the contrary. Since   $A(u)$ form a commuting family, there are an orthogonal $n\times n$-matrix $B$ and a family of
diagonal non-singular real matrices $A'(u)=\mathrm{diag}(a_1(u),\dots,a_n(u))$ such that
$A'(u) = B^{-1}A(u)B$ for $u\in \Omega.$ Consider the bounded and invertible operator
 $(\widehat{B}f)(x):=f(Bx)$ on $L^p(\mathbb{R}^n).$ Because of the equality
 $$
 \widehat{B}\mathcal{H}_{K, A}\widehat{B}^{-1}=\mathcal{H}_{K, A'},
 $$
 operator $\mathcal{H}_{K, A'}$ is Riesz and nontrivial, too.
  For the proof of the Riesz property it is sufficient to verify  the Ruston condition \cite{Ruston} (see also \cite[Theorem 3.12]{Dow}) for this operator. But if $\mathcal{C}$ denotes the ideal of compact operators on $L^p(\mathbb{R}^n)$ one have
 $$
 \inf_{C'\in \mathcal{C}}\|\mathcal{H}_{K, A'}^n-C'\|^{1/n}=\inf_{C'\in \mathcal{C}}\|\widehat{B}\mathcal{H}_{K, A}^n\widehat{B}^{-1}-C'\|^{1/n}
 $$
$$
\le\|\widehat{B}\|^{1/n}\inf_{C'\in \mathcal{C}}\|\mathcal{H}_{K, A}^n-\widehat{B}^{-1}C'\widehat{B}\|^{1/n}\|\widehat{B}\|^{-1/n}
$$
$$
=\|\widehat{B}\|^{1/n}\inf_{C\in \mathcal{C}}\|\mathcal{H}_{K, A}^n-C\|^{1/n}\|\widehat{B}\|^{-1/n}\to 0 \mbox{ as } n\to\infty
$$
which proves the Ruston condition for $\mathcal{H}_{K, A'}$.

As in \cite{Forum} let's consider some  fixed enumeration $U_j\ (j = 1; \dots ; 2^n)$   of the family of all open hyperoctants in $\mathbb{R}^n.$  For every
pair $(i;j)$ of indices there is a unique sequence  $\varepsilon(i; j)\in\{-1,1\}^n$ such that
$$
\varepsilon(i; j)U_i :=
\{(\varepsilon(i; j)_1 x_1;\dots ; \varepsilon(i; j)_n x_n) : x=(x_k)_{k=1}^n\in U_i\} = U_j.
$$
 Then $\varepsilon(i; j)U_j = U_i$ and
$\varepsilon(i; j)U_l \cap U_i=\emptyset$ as $l \ne j.$

 We put
$$
\Omega_{ij} := \{u\in\Omega: (\mathrm{sgn}(a_1(u)); \dots ; \mathrm{sgn}(a_n(u))) = \varepsilon(i; j)\}
 $$
and let
$$
(H_{ij}f)(x) :=\int_{\Omega_{ij}} K(u)f(A'(u)x)d\mu(u).
$$
Since $f(A'(u)x) = 0$ for $f \in L^p(U_j)$
and $x \notin U_i$, each $H_{ij}$ maps $L^p(U_j)$ into $L^p(U_i)$.
Moreover, if $f\in L^p(\mathbb{R}^n)$  and $f_j := f\chi_{U_j}$ ($\chi_E$ denotes the indicator of a subset $E\subset \mathbb{R}^n$) then as in the proof of formula (1) in  \cite{Forum} it is easy to verify that for a.~e. $x \in \mathbb{R}^n$
$$
(\mathcal{H}_{K, A'}f)(x) =\sum\limits_{j=1}^{2^n}\sum\limits_{i=1}^{2^n}(H_{ij}f_j)(x).
$$
  It follows that  for some pare $(i,j)$ of indices the restriction $\mathcal{H}:=H_{ij}|L^p(U_j)$ is a nontrivial operator. In the sequel the pare $(i,j)$ will be fixed.

Consider the map
$Jx:=(\varepsilon(i,j)_kx_k)_{k=1}^n$ $(x=(x_k)_{k=1}^n\in \mathbb{R}^n).$
 Then $J:U_j\to U_i$ and the operator
$(\widehat{J}f)(x):=f(Jx)$
maps $L^p(U_i)$ on $L^p(U_j)$ isometrically. It follows that the operator
$$
\mathcal{K}:=\widehat{J}\mathcal{H}
$$
acts on $L^p(U_j)$ and is bounded. Hereafter we put $U:=U_j$ for simplicity.
   Then $\mathcal{K}$ is a nontrivial Riesz operator on  $L^p(U)$. Indeed, the operator $\widehat{J}\mathcal{H}_{K,A'}=\mathcal{H}_{K,A'}\widehat{J}$ is Riesz on $L^p(\mathbb{R}^n)$ \cite[Lemma 5]{Sche}. Since  the  space $L^p(U)$ is
 invariant with respect to this operator, the restriction $\widehat{J}\mathcal{H}_{K,A'}|L^p(U)$
 (which is equal to  $\mathcal{K}$) is Riesz on  $L^p(U)$  by \cite[p. 80, Theorem 3.21]{Dow}, as well.

Let $1\leq p<\infty.$ To get a contradiction, we shall use the modified $n$-dimensional  Mellin transform $\mathcal{M}$ for the $n$-hyperoctant $U$.
By Lemma 4 the map $\mathcal{M}$ is a bounded operator between  $L^p(U)$ and $L^q(\mathbb{R}^n)$ for $1\leq p\leq 2$  ($1/p+1/q=1).$

  Let $f\in L^p(U)$. Note that
 $$
(\mathcal{K}f)(x) =\int_{\Omega_{ij}} K(u)f(A'(u)(Jx))d\mu(u)=\int_{\Omega_{ij}} K(u)f(A^{\prime\prime}(u)x)d\mu(u)
$$
where
$$
A^{\prime\prime}(u)=\mathrm{diag}(\varepsilon(i,j)_1a_1(u),\dots, \varepsilon(i,j)_na_n(u)).
$$
First assume that   $|y|^{-1/q}f(y)\in L^1(U).$ Then,
  as in the proof of  Theorem 1 from \cite{faa} (or \cite{Forum}), using Fubini--Tonelli's theorem and integrating by substitution     $x=(A(u)^{\prime\prime})^{-1}y$ yield
  $$
  (\mathcal{MK}f)(s)=
\varphi(s)(\mathcal{M}f)(s)\  (s\in \mathbb{R}^n),
$$
where the function
 $$
   \varphi(s):=\int_{\Omega_{ij}}K(u)|a(u)|^{-1/p- i s}d\mu(u)
   $$
(the  $(i;j)$ entry of the matrix symbol of a  Hausdorff operator \cite[Definition 2]{faa}, \cite{Forum}) is bounded and continuous on    $\mathbb{R}^n.$

Thus,
$$
\mathcal{MK}f=\varphi \mathcal{M}f. \eqno(2)
$$
By continuity, the last equality is valid for all  $f\in L^p(U).$

  Let $1\leq p\leq 2.$ There exists a constant $c> 0$ such that the set    $\{s\in \mathbb{R}^n: |\varphi(s)|>c\}$ contains an open ball  $S.$ Formula  (2) implies that
$$
M_{\psi}R_{q,S}\mathcal{M}\mathcal{K}=R_{q,S}\mathcal{M},
$$
  where $\psi=(1/\varphi)|S,$  $M_{\psi}$ denotes the operator of multiplication by  $\psi$, and $R_{q,S}: L^q(\mathbb{R}^n)\to L^q(S),$
  $f\mapsto f|S$ is the restriction  operator. Let $T=R_{q,S}\mathcal{M}.$ Passing to the conjugates  gives
 $$
 \mathcal{K}^*T^*M_{\psi}^*=T^*.
 $$
    By \cite[Theorem 1]{FNR3} this implies that the operator $T^*=\mathcal{M}^*R_{q,S}^*$  has finite rank.
    But by Lemma 3, $R_{q,S}^*$ is the operator of natural embedding $L^p(S)\hookrightarrow L^p(\mathbb{R}^n)$. Thus, the restriction of the operator  $\mathcal{M}^*$  to  $L^p(S)$ has finite rank. Since by Lemma 4 $\mathcal{M}^*$
  can easily be reduced to the Fourier transform, this contradicts the Paley--Wiener theorem on the Fourier image of the space
  $L^2(S)$, see, e.~g., \cite[Theorem  III.4.9]{SW}  (in our  case $L^2(S)\subset L^p(S)$).

  Finally, if $2<p\le \infty$, one can use duality arguments. Indeed, by Lemma 2 the adjoint operator $\mathcal{H}_{K, A'}^*$
  (as an operator on $L^q(\mathbb{R}^n)$) is also of Hausdorff type. More precisely, it equals $\mathcal{H}_{\Psi, B},$ where
  $$B(u)= A(u)'^{-1}= \mathrm{diag}(1/a_1(u),\dots,1/a_n(u))$$  and  $$\Psi(u)=K(u)|\det A(u)'^{-1}|=K(u)/\prod_ja_j(u).$$
  It is easy to verify that  $\mathcal{H}_{\Psi, B}$ satisfies all the  conditions of Theorem 1  (with  $q,$ $\Psi$ and $B$ in place of $p,$ $K$ and $A$, respectively).
   Since $1\le q<2$, the operator  $\mathcal{H}_{\Psi, B}$ is not a Riesz operator on $L^q(\mathbb{R}^n)$. The same is true for
   $\mathcal{H}_{K, A}$,  because   $T$ is a Riesz operator if only if its conjugate $T^*$ is a Riesz operator  \cite[p. 81, Theorem 3.22]{Dow}. This completes the proof.

\section{Corollaries and examples}

\

\textbf{Corollary 1}. \textit{Under the  same assumptions of Theorem 1 the operator $\mathcal{H}_{K,A}-\lambda$ is not Fredholm for some scalar $\lambda \ne 0$.}

\textbf{Corollary 2} \textit{Let the conditions of theorem 1 hold. Then either there is a  non-zero point of $\sigma(\mathcal{H}_{K,A})$ that is not a pole of the resolvent of $\mathcal{H}_{K,A}$, or there is a  non-zero point $\lambda$ of $\sigma(\mathcal{H}_{K,A})$ such that  the spectral projection $E(\lambda)$ has infinite-dimensional range.}

Proof. This follows from Theorem 1 and the characterization of Riesz operator given in \cite[Theorem 3.17]{Dow}.

\textbf{Corollary 3} \cite{Conf}. \textit{Let $A(v)$ be a commuting family of real positive or negative definite $n\times n$-matrices  ($v$ runs over the support of
$K$), and $\det A(v)^{-\frac{1}{p}}K(v)$ $\in L^1(\Omega).$ Then a Hausdorff  operator $\mathcal{H}_{K, A}$
in  $L^p(\mathbb{R}^n)$ ($1\leq p\le\infty$) is non-Riesz  (and, in particular, it is not the sum of a quasinilpotent and compact operators)  provided it is  non-zero.}

\textbf{Corollary 4.} \textit{Let  $\phi:\Omega\to \mathbb{C}$, and  $a:\Omega\to \mathbb{R}$ be  measurable functions, such that $|a(u)|^{-1/p}\phi(u)\in L^1(\Omega).$ Then a one-dimensional  Hausdorff  operator
$$
(\mathcal{H}_{\phi, a}f)(x)=\int_{\Omega} \phi(u)f(a(u)x)d\mu(u)\  (x\in \mathbb{R})
$$
on  $L^p(\mathbb{R})$ ($1\leq p\le\infty$) is non-Riesz   (and, in particular, it is not the sum of a quasinilpotent and compact operators) provided it is  non-zero.}

\textbf{Example 1}. Let $t^{-1/q}\psi(t)\in L^1(0,\infty).$
Then, by Corollary 4, the operator
$$
(\mathcal{H}_\psi f)(x)=\int_0^\infty\frac{\psi(t)}{t}f\left(\frac{x}{t}\right)dt
$$
is a non-Riesz  operator in  $L^p(\mathbb{R})$ ($1\leq p\le\infty$) provided it is non-zero.

\textbf{Example 2}. Let $(t_1t_2)^{-1/p}\psi_2(t_1,t_2)\in L^1(\mathbb{R}_+^2).$
Then, by Corollary 3, the operator
$$
(\mathcal{H}_{\psi_2} f)(x_1,x_2)=\frac{1}{x_1x_2}\int_0^\infty\!\int_0^\infty\psi_2\left(\frac{t_1}{x_1}, \frac{t_2}{x_2}\right) f(t_1,t_2)dt_1dt_2
$$
is a non-Riesz  operator in  $L^p(\mathbb{R}_+^2)$ ($1\leq p\le\infty$) provided it is non-zero.

\textbf{Example 3}. (Discrete   Hausdorff operators, cf.  \cite[Example 3]{Forum}). Let $\Omega=\mathbb{Z}_+^m,$ and $\mu$ be a counting measure.
Then Definition 1 turns into ($f\in L^p(\mathbb{R}^n)$, $1\le p\le\infty$)
$$
(\mathcal{H}_{K, A}f)(x) =\sum_{u\in \mathbb{Z}_+^m} K(u)f(A(u)x)
$$
($A(u)$ form a  family of real non-singular $n\times n$ matrices).
Assume that  $\sum_{u\in \mathbb{Z}_+^m} |K(u)||\det A(u)|^{-1/p}<\infty.$ Then the operator $\mathcal{H}_{K, A}$ is well defined and bounded on
$L^p(\mathbb{R}^n)$ by Lemma 1, and is a non-Riesz  operator   by Theorem 1  provided it is non-zero and matrices $A(u)$ are permutable   and  self-adjoint.

\textbf{Acknowledgment.} This work was supported by the  State Program of Scientific  Research of Republic of Belarus.

 This is a preprint of the article \cite{RJMP}.

\begin{flushleft}

\end{flushleft}

\end{document}